# Learning and Teaching Calculus Through Its History

Chamila Gamage

**Introduction:**

As a calculus instructor at a university, I have always been fascinated by how this subject became one of the most transformative achievements in mathematics. When we look beyond the well-known story of Newton and Leibniz developing calculus independently in the seventeenth century, a much deeper and richer narrative emerges, one that stretches across centuries and cultures and continues to shape how we understand the concepts of change, motion, and growth. Calculus did not appear suddenly or in isolation; it grew out of many traditions in geometry, measurement, and philosophy that together formed the modern language of science.

I see calculus as a truly global achievement rather than a purely European invention. Ideas from ancient Indian, Greek, and Islamic mathematics, along with later European innovations, show how concepts of infinity, area, and continuous change evolved through collaboration and cultural exchange. This historical perspective is valuable not only for instructors but also for students learning mathematics today. When we revisit these challenges, early mathematicians faced in describing infinitesimal change or curved motion, we gain a deeper appreciation for the creativity and persistence behind the formulas we now use so routinely. By linking these historical developments to my classroom practice, I try to present calculus not just as a set of procedures, but as a living story of human thought, and a continuing journey of discovery.

**Cultural Roots of Calculus:**

The history of calculus reaches far beyond the seventeenth century. Its ideas grew slowly from centuries of work by mathematicians in different parts of the world who were all trying to understand motion, area, and change. Long before Newton and Leibniz, ancient scholars in India, Greece, and the Islamic world developed methods that carried the spirit of calculus, even if they lacked its symbols or formal definitions.

In India, the *Śulba-sūtras* (ca. 800–200 BCE), which is considered as the oldest extant Indian mathematical text, explains how geometry was used for practical purposes such as for constructing fire altars with precise shapes and equal areas. As Seshadri 2010 mentions the statement "*the diagonal of a rectangle makes as much (area) as the (areas) made separately by the base and the side put together*" which is a verbal form of the Pythagorean theorem. Later Indian mathematicians continued this pattern of discovery. The Kerala school, led by Madhava in the fourteenth century, produced power-series expansions for sine, cosine, and arctangent functions centuries before Newton and Leibniz (Katz, 1993). These early results highlight the same idea that lies at the heart of calculus: an infinite process can describe a finite value.

Greek mathematics added another crucial piece through the *method of exhaustion* developed by Eudoxus and perfected by Archimedes. By repeatedly inscribing and circumscribing polygons, Archimedes could approximate areas and volumes with any desired accuracy (Heath, 2002). His reasoning shows a deep understanding of limits, even though he expressed it geometrically. This could be a nice way to introduce Riemann sums in a Calculus class, discussing Archimedes' argument for finding the area of a circle or the volume of a sphere and then connect it to the lesson help students see that our modern techniques come from timeless questions about measurement and approximation.

Islamic scholars further advanced these ideas during the medieval period (Rashed, 1994). Ibn al-Haytham, in the eleventh century, derived formulas for sums of powers of integers, relations that later became part of integral calculus. His approach to analyzing motion and light also linked algebra with geometry, showing that concepts like calculus were emerging in multiple regions at once.

In Europe, the fourteenth century French scholar Nicole Oresme brought a new way to visualize these relationships. As Babb (2005) describes, Oresme drew geometric diagrams with one axis representing time and the other representing velocity, creating a visual model that anticipated the modern idea of a function graph. His *mean-speed rule* showed that the distance traveled under constant acceleration equals the area under the velocity-time curve, essentially the first clear statement of what we now call integration. Oresme also studied infinite series, proving that the harmonic series diverges and using geometric figures to represent partial sums.

Together, these examples show that the roots of calculus are truly global. Their combined efforts laid the groundwork for the symbolic calculus that Newton and Leibniz would later formalize.

**The Newton–Leibniz Breakthrough:**

By the seventeenth century, the essential pieces of calculus were already scattered across centuries of mathematical thought. During this period, Isaac Newton in England and Gottfried Wilhelm Leibniz in Germany independently developed what would become known as the calculus (Hathaway, 1919; Katz 1993), a complete system for studying continuous change and accumulation. Their work transformed a diverse collection of ideas into a single language for science and mathematics.

Newton approached calculus through the lens of motion and geometry. As Katz (1993) mentions, his "method of fluxions" described how quantities flow and change over time, much like moving bodies in physics. For Newton, the derivative represented an *instantaneous rate of change*, and the integral measured the *total effect of change* over an interval. These concepts grew naturally from his studies of planetary motion and gravity, later presented in Newton, et al., (1999). Newton's reasoning was highly geometric, expressed through ratios and limits rather than symbols, but it offered an extraordinary new way to connect mathematics with the physical world.

Leibniz, working independently, came to the same discoveries from a very different starting point. His goal was to create a general symbolic language for reasoning about change (Jesseph, 1998). He introduced the notations $dy/dx$ and $\int$, which have since become universal. Leibniz believed that algebraic symbols could capture motion and transformation more precisely than words or diagrams. His notation made calculus teachable and extensible, allowing others to compute systematically (Guicciardini, 1989).

This philosophical viewpoint shaped how calculus developed and was taught. Critics such as Bernard Nieuwentijt and Michel Rolle argued that infinitesimals lacked meaning, since quantities smaller than any real number could not exist (Jesseph, 1998). But Leibniz countered by treating infinitesimals as incomparably small finite quantities, preserving their practical use while avoiding arguments about whether they actually exist. His followers Johann Bernoulli, L'Hôpital, and Varignon adopted this flexible view and turned it into a working method. L'Hôpital's *Analyse des infiniment petits* (1696), the first calculus textbook, presented differentiation rules, examples, and applications in a clear format that defined how calculus would be taught for generations (Blanco, 2024).

The dispute between Newton and Leibniz over who invented calculus became one of the most famous controversies in mathematics. Each claimed credit, supported by national pride and different ways of thinking about the subject (Hathaway, 1919). Yet this had an unexpected benefit, it helped calculus spread quickly across Europe. In England, Newton's geometric method of *fluxions* shaped mathematics for many years, on the other hand, Leibniz's symbolic form became the common language. His notation was short, flexible, and easy to teach, which made it possible for mathematicians to apply calculus to new problems in mechanics, probability, and differential equations (Guicciardini, 1989).

From a teaching perspective, these two approaches, Newton's geometric intuition and Leibniz's symbolic formalism still shape how we teach calculus today. Newton's view speaks to the visual and physical intuition students need when thinking about slopes, velocity, and area. Leibniz's system provides the symbolic precision required for computation and abstraction. Both Newton and Leibniz combined the geometric, algebraic, and philosophical ideas that had been developing for centuries. Their work united ancient intuition with modern reasoning, creating a single way to describe a universe in constant motion.

**Modern Contributions**

The work of Newton and Leibniz gave calculus its form, but their methods were built on intuition. Over the next three centuries, mathematicians refined that intuition into a structure of complete logical rigor. This process transformed calculus from a tool for studying motion into the universal language of modern science.

Leonhard Euler played a decisive role in this evolution. Expanding on Leibniz's notation, Euler unified differentiation and integration through power series and introduced exponential and trigonometric functions in their modern forms. (Katz, 1993)

Joseph-Louis Lagrange attempted to rebuild calculus on purely algebraic foundations. He believed that any differentiable function could be expressed as a Taylor series and that differentiation could be viewed as taking the coefficient of the linear term in that expansion. Though his idea was not entirely successful, it clarified the relationship between derivatives and local linear approximation. Lagrange's focus on structure influenced how calculus would later be presented in classrooms as a systematic theory rather than an art of clever tricks.

The nineteenth century brought the full formalization of calculus through the work of Augustin-Louis Cauchy and Karl Weierstrass. (Katz, 1998) Cauchy introduced precise definitions of limit, continuity, and convergence, and Weierstrass refined them with his ε–δ formulation. Their work gave calculus the rigor it had long lacked, eliminating reliance on infinitesimals and grounding the subject in arithmetic logic. As Zuccheri and Zudini (2014) explain, this shift from geometric intuition to analytical precision fundamentally changed how calculus was taught and understood.

Bernhard Riemann (1854) defined the definite integral as a limit of sums, allowing calculus to handle irregular functions, while Henri Lebesgue (1902) extended this idea by defining integration through measure theory. These contributions opened the door to modern probability, real analysis, and mathematical physics. In parallel, Sophus Lie developed continuous transformation groups that linked calculus with geometry and the emerging theory of symmetry, ideas now essential in physics. David Hilbert and Émile Borel expanded the theoretical foundations of analysis, connecting calculus to topology and infinite-dimensional spaces.

The twentieth century carried calculus into new domains. Norbert Wiener formalized stochastic processes and founded cybernetics, applying calculus to random motion. Alan Turing and others used differential equations to study computation and biological pattern formation. In physics, Einstein's general relativity relied on tensor calculus, while quantum mechanics depended on the differential operators of Schrödinger and Heisenberg. These modern developments show how the spirit of calculus, reasoning about infinitesimal change and continuous transformation remains the core of scientific thought.

For me as an instructor, tracing this journey from intuitive beginnings to modern abstraction has been very inspiring. Each stage in the history of calculus reflects a balance between creativity and discipline, between visual understanding and symbolic reasoning (Ortega, 2002). This same balance defines how students learn. When I teach them from sketching tangent lines to writing limit proofs, or from exploring numerical data to modeling more abstract systems, I see them going through the same path that mathematics itself has taken. The story of calculus is, in many ways, the story of learning itself.

**Integrating Historical Insights into My Teaching**

Studying the history of calculus has made a huge difference in my teaching philosophy. I no longer view my course simply as a sequence of topics like limits, derivatives, integrals, and series, but as a narrative of human discovery that continues in every classroom. By connecting historical developments to modern practice, I could help students see calculus not as a mere collection of facts but as a living subject shaped by curiosity, creativity, and problem solving.

Below are few things I would like to incorporate into my teaching:

1. **Beginning with historical problems.**

For instance, when teaching about area under graphs, we can mention Archimedes' *method of exhaustion*, show them how early mathematicians were motivated, and apply that in the classroom. This approach makes the learning process feel natural as it arises from a real need.

2. **Blending intuition and rigor.**

Following Leibniz's example of "well-grounded fictions," We can encourage students to reason intuitively before introducing formal proofs. At first, we speak of $dy$ and $dx$ as very small changes; later, we refine these ideas through the ε–δ definition of a limit. This mirrors the historical transition from informal reasoning to modern analysis. Students can learn that approximation and understanding often come before full rigor.

3. **Highlighting cross-cultural connections.**

We can make a point of mentioning contributions from Indian, Greek, and Islamic scholars alongside those of Newton and Leibniz. When students see that ideas such as infinite series or geometric algorithms existed long before seventeenth-century Europe, they recognize mathematics as a global conversation. This perspective not only broadens their appreciation of the subject but also invites them to see themselves as potential contributors to those ongoing conversations.

**4. Encouraging reflection.**

We can ask students to write short reflections that connect a historical idea to a modern problem. These reflections encourage them to think about how early mathematicians might have approached a similar challenge, what obstacles they faced, and how those ideas evolved into the methods we use today.

**Conclusion**

Reflecting on the history of calculus has changed how I view both the subject and my role as an instructor. The development of calculus was not a single discovery but the result of many cultures and centuries of thought. Each was motivated by a common desire to understand motion and measure the world more precisely. Their collective effort created a breakthrough that transformed mathematics, connecting geometry, algebra, and analysis into a single framework. Calculus became the foundation for modern science, technology, and the way we model growth, motion, and uncertainty. Learning this history helps me teach calculus as more than a set of procedures. It reminds students that the subject was built through creativity, collaboration, and persistence and that the same qualities are essential in learning it today

(Eds.), *Handbook on the history of mathematics education* (pp. 493–513). Springer. https://doi.org/10.1007/978-1-4614-9155-2_24